\documentclass{amsart}
\usepackage{amsmath}
\usepackage{graphicx}
\usepackage{sidecap}
\usepackage{subfig}

\newtheorem{theorem}{Theorem}
\newtheorem{lemma}{Lemma}
\newtheorem{remark}{Remark}
\newtheorem{definition}{Definition}

\newtheorem{corollary}{Corollary}

\author{Avetik Arakelyan}
\email{arakelyanavetik@gmail.com}
\address{Institute of mathematics,
NAS of Armenia,
Bagramian ave. 24B, 0019 Yerevan, Armenia }

\title [Finite Difference scheme for parabolic Two-phase problem]{A Finite Difference Method for Two-phase  Parabolic Obstacle-like Problem}
\thanks{This work was supported by State Committee of Science MES RA, in frame of the research project No.  SCS 13YR-1A0038}

\begin{document}

\begin{abstract}
In this paper we treat the numerical approximation of the two-phase parabolic obstacle-like problem:
\[\Delta u -u_t=\lambda^+\cdot\chi_{\{u>0\}}-\lambda^-\cdot\chi_{\{u<0\}},\quad (t,x)\in (0,T)\times\Omega,\]
where $T < \infty, \lambda^+ ,\lambda^- > 0$ are Lipschitz  continuous functions, and $\Omega\subset\mathbb{R}^n$ is a bounded domain. We introduce  a certain variational form,  which allows us to define a notion of viscosity solution. We use defined viscosity solutions framework to apply Barles-Souganidis theory. The numerical projected Gauss-Seidel method is constructed. Although the paper is devoted to the parabolic version of the two-phase obstacle-like problem, we prove convergence of the discretized scheme to the unique viscosity solution for both  two-phase parabolic  obstacle-like and  standard two-phase membrane problem. Numerical simulations are also presented.

\end{abstract}

\maketitle

\section{Introduction}
\subsection{The statement and known results}
 In this paper we construct and implement a numerical method for the two-phase parabolic obstacle-like problem
\begin{equation}\label{two_phase_parabolic}
\begin{cases}
\Delta u -u_t=\lambda^+\cdot\chi_{\{u>0\}}-\lambda^-\cdot\chi_{\{u<0\}}, \text{ in}\; (0,T)\times\Omega,\\
u(0,x)=g(x),\qquad\qquad\qquad\qquad\qquad\text{ in}\; \{0\}\times\Omega,\\
u(t,x)=h(t,x),\qquad\qquad\qquad\qquad\qquad\text{ in}\; (0,T)\times\partial\Omega,\\
\end{cases}
\end{equation}
where
\begin{equation}\label{conditions}
  \lambda^\pm\in C^{0,1}(\Omega)\; 0<T<\infty,\;\text{and}\; \Omega\subset \mathbb{R}^n\;\text{is a given bounded
domain.}
\end{equation}

Here $g(x)$ is a sign changing continuous function, and the function $h(t,x)$ is possibly a sign changing function.

 The problem arises as limiting case in the model of temperature control through the interior described in \cite[Section $2.3.2$]{MR0521262}.

 In the paper \cite{MR2511041}  the authors proved an optimal regularity result for the free boundary ${\partial\{u>0\}\cup\partial\{u<0\}}.$ They show that if a branch point occurs (i.e., the two phases $\{u>0\}$ and $\{u<0\}$ coexist with vanishing $\nabla u$), then nearby $\partial\{u>0\}$ and $\partial\{u<0\}$ the boundary is the union of two Lipschitz graphs that are continuously differentiable in the space variable.

The stationary case -\textit{ The two-phase membrane problem}  has been studied from different view points. The optimal $C^{1,1}_{loc}$ regularity has been proved by Ural'tseva \cite{MR1906034} in the case when the coefficients $\lambda^\pm$ are assumed to be constant, and the result was extended by Shahgholian \cite{MR1934623} for Lipschitz-regular  $\lambda^\pm.$ The regularity for the free boundary has been studied by Shahgholian, Ural'tseva and Weiss \cite{MR2258264}, \cite{MR2340105}.
\subsection{Barles-Souganidis theory}
For the sake of clarity in this section we present a very fundamental theorem related to the convergence of monotone difference schemes. The result has been obtained by G.Barles and P.Souganidis in $1991$ (see \cite{MR1115933}) and since then many applications in numerical analysis of \emph{monotone} difference schemes, has been made.

We consider the equations of the following form:
\begin{equation}\label{general_form}
F(D^2u,Du,u_t,u,t,x)=0\quad\text{in}\; [0,T]\times\overline\Omega.
\end{equation}
Here $\Omega$ is an open subset of $\mathbb{R}^n,$ and $\overline\Omega$ is its closure. The functions
\[F:\mathbb{S}^n\times\mathbb{R}^n\times\mathbb{R}\times[0,T]\times\overline\Omega\rightarrow\mathbb{R}\;\;\text{and}\;\;
u:[0,T]\times\overline\Omega\rightarrow\mathbb{R},\]
are bounded (possibly discontiuous), and finally, $Du$ and $D^2u$ stand for the gradient vector and second derivative (Hessian) matrix of $u.$
We say that \eqref{general_form} is \emph{elliptic} if for all $x\in\mathbb{R}^n,t>0,p_t\in\mathbb{R},p_x\in\mathbb{R}^n,\text{and}\quad M,N\in\mathbb{S}^n$
\[F(M,p_x,p_t,u,t,x)\leq F(N,p_x,p_t,u,t,x)\quad\text{ if}\quad M\geq N.\]

Before stating the theorem we need to define some notions related to the finite difference schemes.

A numerical scheme is an equation of the following form
\begin{equation}\label{introd_scheme}
S(h,t,x,u_h(t,x), [u_h]_{t,x}) = 0,
\end{equation}
%for (t; x) in Ghnft = 0g
%uh(0; x) = uh;0(x) in Gh \ ft = 0g
where $u_h$ stands for the approximation of $u$ and $[u_h]_{t,x}$ represents the value of $u_h$ at other points
than $(t,x).$ Here for simplicity we take $\Delta x =\Delta t=h.$
The theory requires the following assumptions:

\textbf{Monotonicity:} If $u\leq v,$
\[S(h, t, x, r, u)\geq S(h, t, x, r, v).\]

\textbf{Consistency:}
For every smooth function $\phi(t,x),$
\[S(h,t,x,\phi(t,x),[\phi(t,x)]_{t,x})\rightarrow F(D^2\phi,D\phi,\phi_t,\phi,t,x),\]
as $\Delta x\rightarrow 0\;\text{and}\; \Delta t\rightarrow 0.$\\

\textbf{Stability:}
For every $h > 0,$ the scheme has a solution $u_h$ which is uniformly bounded independently of $h.$

The theorem reads as follows:

\begin{theorem}(\textbf{Barles-Souganidis 1991})
Under the above assumptions, if the scheme \eqref{introd_scheme} satisfy the \emph{consistency},
\emph{monotonicity} and \emph{stability} property, then its solution $u_h$ converges locally uniformly
to the unique viscosity solution of \eqref{general_form}.
\end{theorem}

\subsection{The outline of the paper}
The paper is organized as follows: In section $2$ we introduce a variational form  which allows us to define a notion of viscosity solutions for the underlying problem.  In section $3$ we  construct a numerical monotone difference scheme using developed varioational form.  Convergence of the scheme to the unique viscosity solution follows from the so-called Barles-Souganidis theorem. In Section $4$ we develop projected Gauss-Seidel algorithm to approximate the discrete two-phase parabolic obstacle-like problem, and  give some numerical examples of the discrete solutions by this algorithm.

\section{Weak and viscosity solutions}
%\subsection{Viscosity solution}
We start this section with definition of viscosity solutions for parabolic type equations:
\begin{equation}\label{GeneralParabol}
G(D^2u,Du,u_t,u,t,x)=0 \quad \text{on}\quad\overline\Omega\times(0,T),
\end{equation}
where $\Omega$ is a bounded domain and $G(M,p_x,p_t,u,t,x)$ is a real-valued discontinuous function defined on $\mathbb{S}^n\times\mathbb{R}^n\times\mathbb{R}\times[0,T]\times\overline\Omega,$ where $\mathbb{S}^n$ is the space of
$n\times n$ symmetric matrices. Here $G$ is always assumed to satisfy the following \emph{ellipticity condition }
\[G(M,p_x,p_t,u,t,x)\leq G(N,p_x,p_t,u,t,x)\quad\text{ if}\quad M\geq N,\]
for all $x\in\mathbb{R}^n,t>0,p_t\in\mathbb{R},p_x\in\mathbb{R}^n,\text{and}\quad M,N\in\mathbb{S}^n$
(see \cite{MR1118699}).
\begin{definition}\label{visc_definition}
A bounded uniformly continuous function $u:[0,T]\times\overline\Omega\rightarrow\mathbb{R}$ is called  a viscosity subsolution (resp. supersolution) for \eqref{GeneralParabol}, if for all $\phi\in C^2([0,T]\times\overline\Omega)$ and all $x\in\Omega$ such that $u-\phi$ has a local maximum (respectively minimum) at $(t,x)$, we have
\[G(D^2\phi,D\phi,\phi_t,u,t,x)\leq 0,\]
\qquad\qquad\qquad (\text{respectively}
\[
G(D^2\phi,D\phi,\phi_t,u,t,x)\geq 0).
\]
\end{definition}
The function $u$ is said to be a \emph{viscosity solution} of \eqref{GeneralParabol} if it is both sub- and supersolution of \eqref{GeneralParabol}. For general background  about the theory of viscosity solutions we refer to \cite{MR1118699}.

Now we define the following variational form:
\begin{equation}\label{variation_form}
\mathcal G[u]\equiv G[D^2u,u_t,u]=\min(u_t-\Delta u +\lambda^+, \max(u_t-\Delta u-\lambda^-,u)),
\end{equation}
where $\lambda^+$ and $\lambda^-$ are positive and Lipschitz continuous as in \eqref{conditions}.
%Now we are ready to define the viscosity solution in our case.
It is easy to see that
\[G(X,p_t,r)=\min(p_t-\text{trace}(X)+\lambda^+,\max(p_t-\text{trace}(X)-\lambda^-,r))\]
satisfies \emph{ellipticity} and other conditions, as stated in the beginning of this section, hence we can apply Definition \eqref{visc_definition} as a notion of viscosity sub- and supersolution.
\begin{lemma}[Uniqueness]\label{unique}
The two-phase parabolic obstacle-like problem \eqref{two_phase_parabolic}  has a unique weak solution.
\end{lemma}
\begin{proof}
Suppose there exist two weak solutions $u(t,x)$ and $v(t,x).$ Then for every $(t,x)\in \{u>v\},$ we have
\[
\chi_{\{u>0\}}\geq\chi_{\{v>0\}}\;\; \mbox{and}\;\; \chi_{\{u<0\}}\leq\chi_{\{v<0\}} .
\]

Thus,
\begin{align*}
\Delta u-u_t &= \lambda^+\cdot\chi_{\{u>0\}}-\lambda^-\cdot\chi_{\{u<0\}}\\&\geq
\lambda^+\cdot\chi_{\{v>0\}}-\lambda^-\cdot\chi_{\{v<0\}}\\&=\Delta v-v_t.
\end{align*}
Therefore  $\Delta(u-v)-(u-v)_t\geq 0$ on the set $\{u>v\}.$  Now, the weak maximum principle clearly gives us $u(t,x)\leq v(t,x),$ which is inconsistent with the set $\{u>v\},$ and hence  $\{u>v\}={\O}.$ Similarly, if we consider the set  $\{u<v\},$ then the same arguments will lead us to   $\{u<v\}={\O}.$ Thus $u(t,x)=v(t,x),$ and this completes the proof of lemma.
\end{proof}

%\begin{lemma}[Comparison principle]\label{comparison}
%Let  $\Omega$ be a bounded domain and $v_1,v_2\in W^{2,\infty}(\Omega)$. If
%$$
% \mathcal G [v_1] \le \mathcal G [v_2]\quad \mbox{a.e. in} \quad \Omega \qquad  \mbox{and} \qquad v_1\le v_2\quad \mbox{on} \quad \partial\Omega,
%$$
%then
%$v_1\le v_2$ in $\Omega.$
%\end{lemma}

%\begin{proof} Let $\Omega_1=\{x\in\Omega:\mbox{  } v_1(x)>v_2(x)\}$.  If the set $\Omega_2=\{x\in \Omega_1: -\Delta v_1(x)>-\Delta v_2(x)\}$ has positive Lebesgue measure, then we get a contradiction, since
%$\mathcal F [v_1](x)> \mathcal F [v_2](x)\quad \mbox{in}\quad \Omega_2.$ Consequently, $-\Delta v_1(x)\leq-\Delta v_2(x)$ a.e. in $\Omega_1.$ But in this case the weak maximum principle implies $v_2\geq v_1$ in $\Omega_1$, which is inconsistent with the definition of $\Omega_1$. Therefore, $\Omega_1=\emptyset$.
%\end{proof}

\begin{theorem}\label{equivalence}
 If $u$ is the solution (in the weak sense) to \eqref{two_phase_parabolic}, then it is a viscosity solution to \eqref{variation_form} (with the same boundary conditions as in \eqref{two_phase_parabolic})  and vice versa.
\end{theorem}
\begin{proof}
Suppose $u^*$ solves the two-phase parabolic equation \eqref{two_phase_parabolic}. Then  \eqref{two_phase_parabolic} will satisfy the following inequality in the sense of distributions
\[-\lambda^-\leq\Delta u^*-u^*_t \leq\lambda^+\quad \text{ in}\;[0,T]\times\Omega,\]
 and hence  it holds in the viscosity sense as well (see \cite{MR1341739}).
We consider two cases:\\
\begin{itemize}
%\item{$(t,x)\in\{u^*>0\}$}\\
%In this case we note that the solution will be a smooth function, in fact $u^*\in C^{2,1}_{x,t}(\{u^*>0\}),$ hence
%one can understand the derivatives in the classical sense. Apparently \eqref{two_phase_parabolic} will be reduced to
%\[u^*_t-\Delta u^* +\lambda^+=0.\]
%On the other hand \[\max(u^*_t-\Delta u^*-\lambda^-,u^*)>0,\] because of $u^*(t,x)>0$. Therefore one obtains \[G[D^2u^*,u^*_t,u^*]=\min(0,\max(u^*_t-\Delta u^*-\lambda^-,u^*))=0.\]
%Hence $u^*$ solves \eqref{variation_form} as well.\\

\item{$(t,x)\in\{u^*=0\}$}\\
In this case as mentioned above the equation \eqref{two_phase_parabolic} will satisfy the following inequality in the viscosity sense
\begin{equation}\label{visc_ineq}
-\lambda^-\leq\Delta u^*-u^*_t \leq\lambda^+,\;\; \text{in}\;[0,T]\times\Omega.
\end{equation}
Let $(t_0,x_0)\in\{u^*=0\}$ and $\phi\in C^2([0,T]\times\Omega),$ such that $u^*-\phi$ has a local minimum at $(t_0,x_0).$
Then  according to \eqref{visc_ineq} and definition of viscosity supersolution we obtain
\[\Delta\phi(t_0,x_0)-\phi_t(t_0,x_0)\leq\lambda^+(t_0,x_0).\]
Since
 \[\max(\phi_t(t_0,x_0)-\Delta \phi(t_0,x_0)-\lambda^-(t_0,x_0),0)\geq 0,\]
then we easily obtain
\[\qquad\qquad G[D^2\phi,\phi_t,u^*](t_0,x_0)\geq 0.\]
Thus $u^*$ is a viscosity supersolution  for our variational equation for all $(t,x)\in\{u^*=0\}.$
In the same way if $(t_0,x_0)\in\{u^*=0\}$ and $\psi\in C^2([0,T]\times\Omega),$ such that $u^*-\psi$ has a local maximum at $(t_0,x_0),$ then by taking again into account \eqref{visc_ineq}, we have
\[\Delta\psi(t_0,x_0)-\psi_t(t_0,x_0)\geq-\lambda^-(t_0,x_0),\]
and hence we obtain
\[G[D^2\psi,\psi_t,u^*](t_0,x_0)\leq 0.\]
Thus $u^*$ is also a viscosity subsolution at that point. Therefore $u^*$ is a viscosity solution to  \eqref{variation_form}, for all $(t,x)\in\{u^*=0\}$.\\
%and hence it holds in viscosity sense as well (see \cite{MR1341739}).
%therefore
%\[\max(u^*_t-\Delta u^*-\lambda^-,0)=0,\]
%holds in the viscosity sense as well.
%Thus
%\begin{align*}
%G[D^2u^*,u^*,u^*_t] &=\min(u^*_t-\Delta u^*+\lambda^+,\max(u^*_t-\Delta u^*-\lambda^-,0))\\
%                    &=\min(u^*_t-\Delta u^*+\lambda^+,0)=0,
%\end{align*}
%in the viscosity sense.

\item {$(t,x)\in\{u^*<0\}\cup \{u^*>0\}$}\\
Note that in this case the solution $u^*$ will be $C^{2,1}_{x,t}$ smooth in a small neighborhood of the point $(t,x).$  Thus, one can understand the derivatives in the classical sense. Now, if we assume, without loss of generality, that $(t,x)\in\{u^*<0\},$ then the equation \eqref{two_phase_parabolic} will be reduced to
\[u^*_t-\Delta u^* -\lambda^-=0.\]
Variational equation \eqref{variation_form} will lead us to
\begin{align*}
G[D^2u^*,u^*_t,u^*] &=\min(u^*_t-\Delta u^* +\lambda^+,\max(0,u^*))\\
                    &=\min(u^*_t-\Delta u^* +\lambda^+,0)\\
                    &=\min(\lambda^-+\lambda^+,0)=0.
\end{align*}
Hence, in this case $u^*$ solves \eqref{variation_form}.
For the case  $(t,x)\in\{u^*>0\},$ we will arrive at
 \[u^*_t-\Delta u^* +\lambda^+=0.\]
On the other hand \[\max(u^*_t-\Delta u^*-\lambda^-,u^*)>0,\] because of $u^*(t,x)>0$. Therefore one obtains \[G[D^2u^*,u^*_t,u^*]=\min(0,\max(u^*_t-\Delta u^*-\lambda^-,u^*))=0.\]
Hence, again  $u^*$ solves \eqref{variation_form}.\\
\end{itemize}

%\end{proof}
%\begin{theorem}[Uniqueness of viscosity solution]
%There exists at most one viscosity solution of   \eqref{variation_form}.
%\end{theorem}
%\begin{proof}
%Suppose there exist two viscosity solutions $u$ and $v$. We have
Next, we are going to proof the opposite statement of the theorem, namely if $u$  is a viscosity solution to
\begin{equation}\label{visc_form_u}
\begin{cases}
\min(u_t-\Delta u +\lambda^+, \max(u_t-\Delta u-\lambda^-,u))=0, \text{ in}\; (0,T)\times\Omega,\\
u(0,x)=g(x),\qquad\qquad\qquad\qquad\qquad\text{ in}\; \{0\}\times\Omega,\\
u(t,x)=h(t,x),\qquad\qquad\qquad\qquad\qquad\text{ in}\; (0,T)\times\partial\Omega,\\
\end{cases}
\end{equation}
 then  it is a weak solution to our two-phase parabolic obstacle-like problem \eqref{two_phase_parabolic}. To this aim,
%\begin{equation}\label{visc_form_v}
%\min(v_t-\Delta v +\lambda^+, \max(v_t-\Delta v-\lambda^-,v))=0,
%\end{equation}
%in the viscosity sense.
%Denote $\tilde\Omega_T=\{(t,x) |\; u(t,x)>v(t,x)\}.$ It is easy to see that
%\begin{equation}\label{Sets_ union}
%\tilde\Omega_T=\{u>v\geq 0\}\cup\{u>0\geq v\}\cup\{u\geq 0>v\}\cup\{0\geq u>v\}.
%\end{equation}
since $u$ satisfies \eqref{visc_form_u}, then  for all $\phi\in C^2(\overline\Omega\times[0,T])$ we will have
\[\min(\phi_t-\Delta \phi +\lambda^+, \max(\phi_t-\Delta \phi-\lambda^-,u))\geq 0,\]
whenever  $u-\phi$ has a local minimum at $(t,x)$. Thus
$\phi_t-\Delta \phi +\lambda^+\geq 0,$ therefore according to definition \eqref{visc_definition}  we conclude that $u$ is a viscosity supersolution to
\[u_t-\Delta u +\lambda^+=0,\]
over $(t,x)\in\overline[0,T]\times\Omega.$ Thus, $u_t-\Delta u +\lambda^+\geq 0$ in the viscosity sense.
On the other hand  \eqref{visc_form_u} implies also
\[\min(\psi_t-\Delta \psi +\lambda^+, \max(\psi_t-\Delta \psi-\lambda^-,u))\leq 0,\]
whenever  $u-\psi$ has a local maximum at $(t,x)$, and $\psi\in C^2(\overline\Omega\times[0,T]).$ This particularly yields
$$
\psi_t-\Delta \psi-\lambda^-\leq\min(\psi_t-\Delta \psi +\lambda^+, \max(\psi_t-\Delta \psi-\lambda^-,u))\leq 0,
$$
which implies that $u$ is also a viscosity subsolution to the equation $u_t-\Delta u -\lambda^-=0.$ The latter statement can be written as $u_t-\Delta u -\lambda^-\leq 0$ in the viscosity sense. Thus we have the following chain of inequalities
\[-\lambda^+\leq u_t-\Delta u\leq\lambda^-,\]
in the viscosity sense, which in turn implies that the inequalities hold in the sense of distributions as well.

Now, by assuming $u(t,x)>0$ we conclude  \[\max(\psi_t-\Delta \psi-\lambda^-,u)>0,\] therefore
$\psi_t-\Delta \psi +\lambda^+ \leq 0.$ This in turn implies that $u$ is a viscosity subsolution to
\[u_t-\Delta u +\lambda^+=0.\]
Hence we conclude that $u$ is a viscosity solution to  $u_t-\Delta u +\lambda^+=0,$
for all $(t,x)\in\{u>0\}.$ Thus
\[u_t-\Delta u=-\lambda^+,\]
in the viscosity sense whenever $u(t,x)>0.$ Recalling again the equivalence property of distributional and viscosity  solutions  for linear PDEs \cite{MR1341739}, we conclude that $u$ is also a distributional solution to  $u_t-\Delta u=-\lambda^+$ on the set $\{u(t,x)>0\}.$

Similarly we will obtain that
\[u_t-\Delta u=\lambda^-,\]
in the viscosity sense, whenever $u(t,x)<0,$ and thus it holds in the distributional sense as well. Thus, we have that the viscosity solution $u$ satisfies the following system in the distributional sense
\begin{align}\label{distrib_cond}
\begin{cases}
\Delta u-u_t=\lambda^+&\mbox{in}\;\;\{u>0\},\\
\Delta u-u_t=-\lambda^- &\mbox{in}\;\;\{u<0\},\\
-\lambda^-\leq \Delta u-u_t\leq\lambda^+&\mbox{ in}\;(0,T)\times\Omega,\\
u(0,x)=g(x),&\mbox{in}\;\{0\}\times\Omega,\\
u(t,x)=h(t,x),&\mbox{in}\;(0,T)\times\partial\Omega.
\end{cases}
\end{align}

Let us assume that $w$ is a unique weak (distributional) solution (Due to Lemma \ref{unique}) to the equation \eqref{two_phase_parabolic}. Then, we have that for every test function $\varphi\in C_0^{\infty}(\overline\Omega\times[0,T])$ the following equality holds
\[
-\int_{\Omega\times(0,T)}\nabla w\cdot\nabla\varphi+\int_{\Omega\times(0,T)}w\cdot\varphi_t=\int_{\Omega\times(0,T)}\left(\lambda^+\cdot\chi_{\{w>0\}}\varphi-\lambda^-\cdot\chi_{\{w<0\}}\varphi\right).
\]
It is clear that for every test function $\varphi\in C_0^{\infty}(\{w>0\})\subset C_0^{\infty}(\overline\Omega\times[0,T])$  we obtain
\[
-\int_{\{w>0\}}\nabla w\cdot\nabla\varphi+\int_{\{w>0\}}w\cdot\varphi_t=\int_{\{w>0\}}\lambda^+\cdot\varphi,
\]
which is nothing else but  the  equation $ \Delta w-w_t=\lambda^+$ in the distributional sense over the set $\{w>0\}.$ Similarly,  if we take a test function $\varphi\in C_0^{\infty}(\{w<0\})\subset C_0^{\infty}(\overline\Omega\times[0,T])$ we will obtain  the  equation $ \Delta w-w_t=-\lambda^-$ in the distributional sense over the set $\{w<0\}.$ On one hand, for every nonnegative test function   $\varphi\in C_0^{\infty}(\overline\Omega\times[0,T])$ we have
\[
-\int_{\Omega\times(0,T)}\nabla w\cdot\nabla\varphi+\int_{\Omega\times(0,T)}w\cdot\varphi_t\geq-\int_{\Omega\times(0,T)}\lambda^-\cdot\chi_{\{w<0\}}\varphi\geq-\int_{\Omega\times(0,T)}\lambda^-\varphi,
\]
and on the other hand for every nonpositive test function  $\varphi\in C_0^{\infty}(\overline\Omega\times[0,T])$  we have
\[
-\int_{\Omega\times(0,T)}\nabla w\cdot\nabla\varphi+\int_{\Omega\times(0,T)}w\cdot\varphi_t\geq\int_{\Omega\times(0,T)}\lambda^+\cdot\chi_{\{w>0\}}\varphi\geq\int_{\Omega\times(0,T)}\lambda^+\varphi,
\]
which leads to the following chain of inequalities $-\lambda^-\leq \Delta w-w_t\leq\lambda^+$ in the sense of distributions. Thus,  the weak solution $w$ satisfies the system \eqref{distrib_cond} as well. Now, if we prove that the system  \eqref{distrib_cond} has a unique solution we are done. We use the same comparison principle as in Lemma \ref{unique}. Let $u$ and $v$ satisfy the system \eqref{distrib_cond}, then if we assume    $\{u>v\}\neq{\O},$  and taking into account that  $$\{u>v\}=\{u>v\geq 0\}\cup\{u>0\geq v\}\cup\{u\geq 0>v\}\cup\{0\geq u>v\},$$
 we obtain the following inequality in the sense of distributions
 $$
 \Delta u-u_t\geq \Delta v-v_t,
 $$
 over the set $\{u>v\}.$ Therefore  $\Delta(u-v)-(u-v)_t\geq 0$ on the set $\{u>v\}.$  Now, the weak maximum principle  gives us $u(t,x)\leq v(t,x),$ which is inconsistent with the set $\{u>v\},$ and hence  $\{u>v\}={\O}.$ Similarly, if we consider the set  $\{u<v\},$ then the same arguments will lead us to   $\{u<v\}={\O}.$ Thus $u(t,x)=v(t,x).$ % \cite[Chapter $3$]{chipot1984}.
\end{proof}

\begin{corollary}
Note that according to  Theorem \ref{equivalence} and Lemma \ref{unique} the viscosity solution $u$ to the variational equation \eqref{visc_form_u} must be unique.
\end{corollary}

\begin{remark}
We note that the same approach works for the two-phase membrane problem as well. In this case the variational equation will be
\begin{equation}\label{Two_phase_membrane}
F[D^2u,u]=\min(-\Delta u +\lambda^+, \max(-\Delta u-\lambda^-,u))=0.
\end{equation}
For more information about the two-phase membrane problem we refer to the papers \cite{MR2340105,MR2258264,MR1906034,MR1825655}. For numerical analysis we refer to the papers \cite{MR1115934,MR2961456,MR1111111}.
\end{remark}
\section{Convergence of numerical schemes}

\subsection{Convergence for the two-phase membrane problem}
$\;\;$Let $\Omega$ be a bounded domain in $\mathbb{R}^n$, $Du$ and $D^2 u$
denote the gradient and Hessian of u, respectively, and $F(x, r, p,X)$ be a continuous
real valued function defined on $\Omega\times\mathbb{R}\times\mathbb{R}^n\times\mathbb{S}^n$, with $\mathbb{S}^n$ being the space of symmetric
$n\times n$ matrices. Write
\[F[u](x)\equiv F(x, u(x),Du(x),D^2u(x)).\]
 Consider the nonlinear, degenerate elliptic partial differential equation with Dirichlet boundary conditions,
\begin{equation}\label{PDE}
\begin{cases}

F[u](x) = 0 \quad   \mbox {in } \Omega,\\
u(x) = g(x)  \quad  \mbox {on } \partial\Omega.
\end{cases}
\end{equation}
  %A. Oberman (see \cite {MR2218974})
We define a uniform  structured grid on the domain $\Omega$ as a directed graph consisting of a set of points $x_i \in \Omega$, $i = 1, . . . ,N,$   each endowed with a number of neighbors $K$.
A grid function is a real valued function defined on the grid, with values $u_i = u(x_i).$
The typical examples of such grid are 3-point and 5-point stencil discretization for the spaces of one dimension and two dimension, respectively.

A function $F^h:\mathbb{R}^N\rightarrow\mathbb{R}^N,$ which is regarded as a map from grid functions to grid functions, is a \emph{finite difference scheme} if
\[F^h[u]^i = F^i[u_i, u_i-u_{i_1},\dots,u_i-u_{i_{K}}] \quad (i = 1,\dots,N),\]
where $\{i_1,i_2,\dots,i_{K}\}$ are the neighbor points of a grid point $i.$ Denote
\[
F^i[u]\equiv F^i[u_i,u_i-u_{i_j}|_{j=\overline{1,K}}]\equiv F^i[u_i,u_i-u_j],\;\;  i = 1, . . . ,N ,
\]
where $u_j$ is shorthand for the list of neighbors $u_{i_j}|_{j=\overline{1,K}}$.

\begin{definition}\label{def_deg_elliptic_scheme}
 The scheme $F$ is degenerate elliptic if each component $F^i$ is nondecreasing
in each variable, i.e. each component of the scheme $F^i$ is a nondecreasing
function of $u_i$ and the differences $u_i - u_{i_j}$ for $j = 1, . . . ,K$.
\end{definition}

Since the grid is uniformly structured, we denote  $h>0$ be the size of the mesh. Then in our case the approximation scheme for the elliptic two-phase membrane
problem is
\begin{equation}\label{numscheme}
F^i[u_i,u_i-u_j]=\min\left(L_h u_i+\lambda^+_i,\max\left(L_h u_i-\lambda^-_i,u_i\right)\right),
\end{equation}
where
\begin{equation}
L_h u_i=\sum_{j=1}^{K}\frac{1}{h^2}(u_i-u_{i_j}),\;\;  i = 1, . . . ,N .
\end{equation}
It is easy to see that $F^i[u_i,u_i-u_j]$ is non-decreasing with respect to $u_i$ and $u_i-u_j$, therefore the finite difference scheme for two-phase membrane problem is a \emph{degenerate elliptic scheme}. But we know that the degenerate elliptic schemes are\emph{ monotone} and \emph{stable} (see \cite {MR2218974}).
\begin{definition}(Consistency).
We say the scheme $F^h$ is consistent with the
equation \eqref{PDE} at $x_0$ if for every twice continuously differentiable function $\phi(x)$
defined in a neighborhood of $x_0$,
\[
F^h[\phi](x_0) \rightarrow F[\phi](x_0) \quad as\quad  h\rightarrow 0.
\]
The global scheme defined on $\Omega$ is consistent if the limit above holds uniformly for all $x\in\Omega.$
(The domain is assumed to be closed and bounded).
\end{definition}
\begin{lemma}\label{consist}
The approximation scheme \eqref{numscheme} is consistent.
\end{lemma}
\begin{proof}
In order to show the consistency we will apply approximation scheme \eqref{numscheme} to the twice continuously differentiable function $\phi(x).$

Suppose $x_0$ is a grid point and the function $\phi(x)$ is twice continuously differentiable around that point.
Then if we use  Taylor expansion for the function $\phi(x)$ around the point $x_0,$ we obtain
\[
L_h \phi_i=\sum_{j=1}^{K}\frac{1}{h^2}(\phi_i-\phi_{i_j})\rightarrow -\Delta\phi(x_0),
\]
and
\[\phi_i\rightarrow\phi(x_0)\quad as\quad h\rightarrow 0. \]
Therefore
\[
\max(L_h \phi_i - \lambda^-_i,\phi_i)\rightarrow\max(-\Delta \phi(x_0) - \lambda^-(x_0),\phi(x_0)),
\]
and
\[
L_h \phi_i + \lambda^+_i\rightarrow-\Delta\phi(x_0) + \lambda^+(x_0)
\]
as $h\rightarrow 0.$
%\begin{equation}\label{consistmin}
%\min(L_h \phi_i + \lambda^+_i ,\max(L_h \phi_i - \lambda^-_i,\phi_i))\rightarrow\min(-\Delta\phi(x_0) + \lambda^+(x_0) ,\max(-\Delta \phi(x_0) - \lambda^-(x_0),\phi(x_0))
%\end{equation}
%as $h\rightarrow 0.$
Thus in light of \eqref{Two_phase_membrane} and \eqref{numscheme}
\[
F^h[\phi](x_0) \rightarrow F[\phi](x_0) \quad as\quad  h\rightarrow 0.
\]
\end{proof}
Now we are ready to formulate the convergence result for the  two-phase membrane problem.
\begin{theorem}(Convergence)
 The finite difference scheme given by \eqref{numscheme} converges uniformly on compacts
subsets of  $\Omega$ to the unique viscosity solution of the two phase-membrane variational equation \eqref{Two_phase_membrane}.
\end{theorem}

\begin{proof}
By virtue of the so-called Barles-Souganidis Theorem (see \cite{MR1115933}) we need  to show that the scheme
is monotone, stable and consistent. The stability and monotonicity are provided by  Definition  \ref{def_deg_elliptic_scheme}, and the finite difference scheme \eqref{numscheme}, where we have concluded that it is actually a degenerate elliptic scheme. Consistency follows from lemma \ref{consist}.
This completes the proof of the theorem.
\end{proof}
\begin{remark}
Note that for the two phase membrane problem the corresponding min-max variational form was introduced in \cite{MR2961456}.
\end{remark}

\subsection{Convergence for Parabolic two-phase obstacle-like problem}
$\;$ Define $\Omega_T=\Omega\times(0,T).$ As in previous section we consider a uniform structured grid on the domain $\Omega$ consisting of a set of points $x_i \in \Omega$, $i = 1, . . . ,N,$  with a number of neighbors $K.$ For the time axis discretization we use the following grid $t_j \in [0,T]$, where  $j = 1, . . . ,M$.

Unfortunately the notion of \textit{degenerate elliptic schemes} is not applicable to this case, since the variational form \eqref{variation_form} is not degenerate parabolic. To be more clear it is noteworthy that  the degenerate elliptic schemes defined above and applied for elliptic version of the Two-Phase Obstacle-like problem,  are just particular case of the schemes defined in the Barles-Souganidis result. But since the variation form of Parabolic two-phase obstacle-like equation is not degenerate parabolic as stated above, then we have to proceed all the steps to check whether our scheme satisfies the required conditions for monotone schemes stated in the Barles-Souganidis theorem or not.
In order to do that we  follow the notations of \cite{MR1315000}. A numerical scheme can be written as
\[S(m,\tilde u)\equiv S(\Delta t, \Delta x, m, j, u_j^m, \tilde u)=0,\]
for $1\leq j\leq N$ and $1\leq m\leq M,$ where $N$ and $M$ are respectively the number of grid points in space and in time.
Here $\tilde u$ denotes the vector $(u_k^l)_{k,l}$ in $\mathbb{R}^{NM}.$ Finally $\Delta t$ and $\Delta x$ denote the time and the space mesh size respectively.

The definition of \textit{Monotonicity} for the scheme will be as follows:
\[
S(\Delta t, \Delta x, m+1, j, u_j^{m+1}, \tilde u)\leq S(\Delta t, \Delta x, m+1, j, v_j^{m+1},
 \tilde v)\quad\text{if}\;\tilde u\geq\tilde v\;\text{and if}\; u_j^{m+1}=v_j^{m+1},
\]
for any $\Delta t,\Delta x>0, 1\leq j\leq N, 1\leq m\leq M,$ and for all $\tilde u$ and $\tilde v$ in  $\mathbb{R}^{NM}.$
In our case the scheme will have the following form
\begin{equation}\label{scheme_parabolic}
S(m+1,\tilde u)\equiv\min(\tilde S(\tilde u)+\Delta t {\lambda^+}_j;
 \max(\tilde S(\tilde u)-\Delta t {\lambda^-}_j, \Delta t u_j^{m+1}) )=0,
\end{equation}
where
\[\tilde S(\tilde u)=u_j^{m+1}-u_j^{m}+\frac{\Delta t}{(\Delta x)^2}Lu_j^m,\; \text{and}\;
Lu_i^m=\sum_{q=1}^{K}(u_i^m-u_{i_q}^m),\;\;  i = 1,\dots,N, \]
here $\{i_1,\dots,i_K\}$ are the neighbor points of a grid point $i.$

\begin{lemma}\label{lemma_parabolic}
 The scheme \eqref{scheme_parabolic} is \textit{monotone} and \textit{stable} provided the following non-linear CFL condition holds
\begin{equation}
\frac{\Delta t}{(\Delta x)^2}\leq\frac{1}{K}.
\end{equation}
\end{lemma}
\begin{proof}
In order to prove monotonicity of \eqref{scheme_parabolic} it is enough to see that \[\tilde S(\tilde u)\leq\tilde S(\tilde v),\quad\text{if}\;\tilde u\geq\tilde v\;\text{and if}\; u_j^{m+1}=v_j^{m+1},\]
for any $\Delta t,\Delta x>0, 1\leq j\leq N, 1\leq m\leq M,$ and for all $\tilde u$ and $\tilde v$ in
$\mathbb{R}^{NM}.$

We have
\[\tilde S(\tilde u)=u_j^{m+1}-u_j^{m}+\frac{\Delta t}{(\Delta x)^2}\sum_{q=1}^{K}(u_j^m-u_{j_q}^m)=u_j^{m+1}-u_j^m\left(1-\frac{\Delta t}{(\Delta x)^2}K\right)-\frac{\Delta t}{(\Delta x)^2}\sum_{q=1}^{K}u_{j_q}^m.\]
Since $u_j^{m+1}=v_j^{m+1}$ and $u_j^{m}\geq v_j^{m}$ for all $1\leq j\leq N, 1\leq m\leq M$ and from CFL condition
\[ K\frac{\Delta t}{(\Delta x)^2}\leq 1,\]
we arrive at
\[\tilde S(\tilde u)\leq\tilde S(\tilde v).\]

To obtain the stability we refer to the Lemma 4.1 in \cite{MR2532350}, where the authors proved comparison principle of numerical scheme defined for one-phase parabolic type equation arising in American option valuation. It is not hard to see that we can do the same induction for our scheme as well (this is standard). Once we have this, the stability follows directly, because we can mimic with the boundary values of the scheme and see that for fixed spatial and time discretization our scheme  will stay between the maximum and minimum values of the discrete boundary values which are fixed a priori. This part is also standard to proceed for such schemes, that's why we skip the detailed proof.

Thus we have the stability and monotonicity for \eqref{scheme_parabolic}.
\end{proof}
\textit{Consistency} of \eqref{scheme_parabolic} can be done as in previous section for the two-phase membrane problem.
We can easily observe that the following limit holds
\begin{equation}\label{consist_parabolic}
\frac{S(\Delta t, \Delta x, m+1, j, \phi_j^{m+1}, \tilde \phi)}{\Delta t}\rightarrow G[D^2\phi,\phi,\phi_t],
\end{equation}
as $\Delta t,\Delta x\rightarrow 0,$ for every twice continuously differentiable function $\phi(t,x).$
\begin{remark}
It is easy to see that  $3-$point and $5-$point stencil discretization will lead us to the following CFL conditions \[\frac{\Delta t}{(\Delta x)^2}\leq\frac{1}{2}\;\;\;\text{and}\;\;\;\frac{\Delta t}{(\Delta x)^2}\leq\frac{1}{4},\;\;\;\text{respectively}.\]
\end{remark}
\begin{remark}
We note that we can consider the implicit discretization of this scheme as well. In this case we will have \emph{unconditionally} monotone and stable scheme.
\end{remark}
We are ready to write down the main result of this section.
\begin{theorem}(Convergence for parabolic case)
The solution $\tilde u$ of\eqref{scheme_parabolic} converges as $\Delta t,\Delta x\rightarrow 0$ uniformly on compacts
subsets of  $\Omega_T$ to the unique viscosity solution of the two-phase parabolic obstacle-like  variation equation \eqref{visc_form_u} .
\end{theorem}
\begin{proof}
The proof is again immediate consequence of Barles-Souganidis theorem as in previous section. The stability and monotonicity are provided by lemma \eqref{lemma_parabolic}. The consistency follows from the limit \eqref{consist_parabolic}.
\end{proof}

\section{Numerical method and Simulations}
\subsection{Numerical method}
As mentioned in the abstract of the paper for constructing a numerical method we refer to Nonlinear Gauss-Seidel method.

Suppose $u^m$ is a shorthand of $(u_j^m)_j.$ We proceed as follows:
\begin{itemize}
\item{\textbf{First Step}.}
\[u^{m+\frac{1}{2}}=\min\left(u^m-\frac{\Delta t}{(\Delta x)^2}Lu^m+\Delta t\lambda^-,0\right),\]
\item{\textbf{Second Step}.}
\[u^{m+1}=\max\left(u^m-\frac{\Delta t}{(\Delta x)^2}Lu^m-\Delta t\lambda^+,u^{m+\frac{1}{2}}\right).\]
\end{itemize}
In order to see the consistency of the method with the  difference scheme \eqref{scheme_parabolic}, we eliminate $u^{m+\frac{1}{2}}$ in the above equality. We obtain
\[u^{m+1}=\max\left(u^m-\frac{\Delta t}{(\Delta x)^2}Lu^m-\Delta t\lambda^+,\min\left(u^m-\frac{\Delta t}{(\Delta x)^2}Lu^m+\Delta t\lambda^-,0\right)\right).\]
Therefore
\[\min\left(u^{m+1}-u^m+\frac{\Delta t}{(\Delta x)^2}Lu^m+\Delta t\lambda^+,\max\left(u^{m+1}-u^m+\frac{\Delta t}{(\Delta x)^2}Lu^m-\Delta t\lambda^-,u^{m+1}\right)\right)=0.\]
Dividing the first argument in the \textit{max} and in the \textit{min} by $\Delta t$ we will derive the desired consistency condition.
\subsection{Simulations}
In this section we present some numerical examples for the two-phase parabolic obstacle-like problem. For all examples we consider $\Omega=[0,1]$ and $T=1.$ %For different $\lambda^\pm$ numerical simulations are shown.

Our equation reads:
\begin{equation}\label{parabolic_example}
\left \{
\begin{array}{ll}
u_{xx}-u_t=\lambda^+\chi_{\{u>0\}} -\lambda^-\chi_{\{u<0\}}, & (t,x)\in (0,1)\times(0,1), \\
u(0,x)=g(x),\\
u(t,0)=h_1(t),u(t,1)=h_2(t).\\
\end{array}
\right.
\end{equation}
\[\]

In the figures \ref{1,3}-\ref{6,6} for different $\lambda^+$ and $\lambda^-$ numerical simulations are shown.
For all cases we take the initial data at time $t=0$ to be linear and the boundary values $h_1(t)$ and $h_2(t)$ to be constant. We use \emph{implicit} discretization in space and forward discretization in time. Numerical examples were constructed with 200 discretization points in space and 250 discretization points in time. In the contour plots  of figures  are clearly visible the positivity and negativity sets of solutions. In all cases we observe that as time evolves, the free boundary becomes more stable and after some amount of time it does not change much. This is expected since from the theory of Parabolic Two-phase Obstacle-Like problems we know that $||u_t||_{\infty}\to 0,$ as $t\to\infty.$

%In the figure \ref{2,7,16} we take $\lambda^+=0.7$ and $\lambda^-=0.2.$ Here our initial data $u(0,x)=16(x-1/2)^3.$ And again the solution was constructed with 200 discretization points in space and 250 discretization points in time. The figure clearly shows that the zero set is bigger than in previous pictures, this is due to the initial data where we have more curvature.

 \begin{figure}[!htbp]
%{\includegraphics[width=7.0cm]{triangel}}
%\centering
 \begin{center}
  \subfloat[Numerical solution]{\includegraphics[width=.48\textwidth]{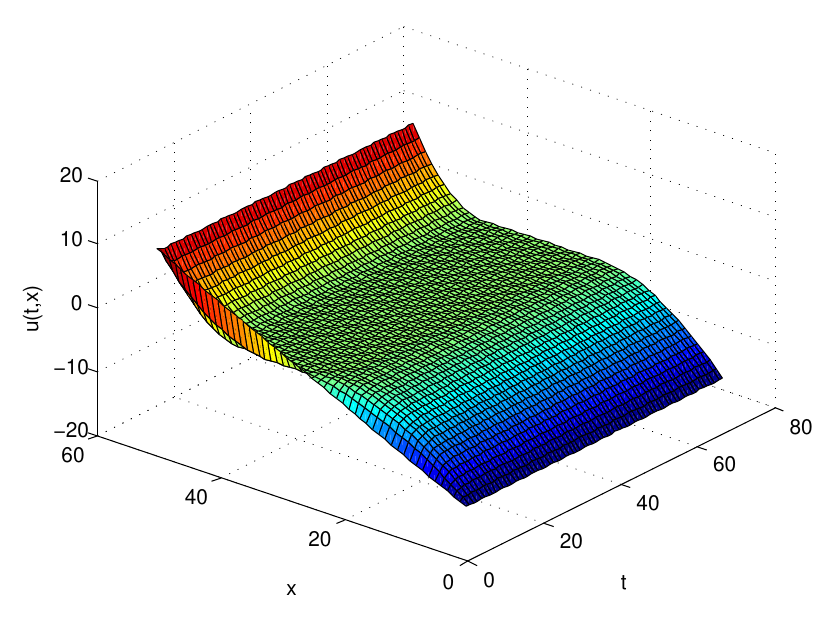}}
  \hspace{1cm}
  \subfloat[Contours]{\includegraphics[width=.4\textwidth]{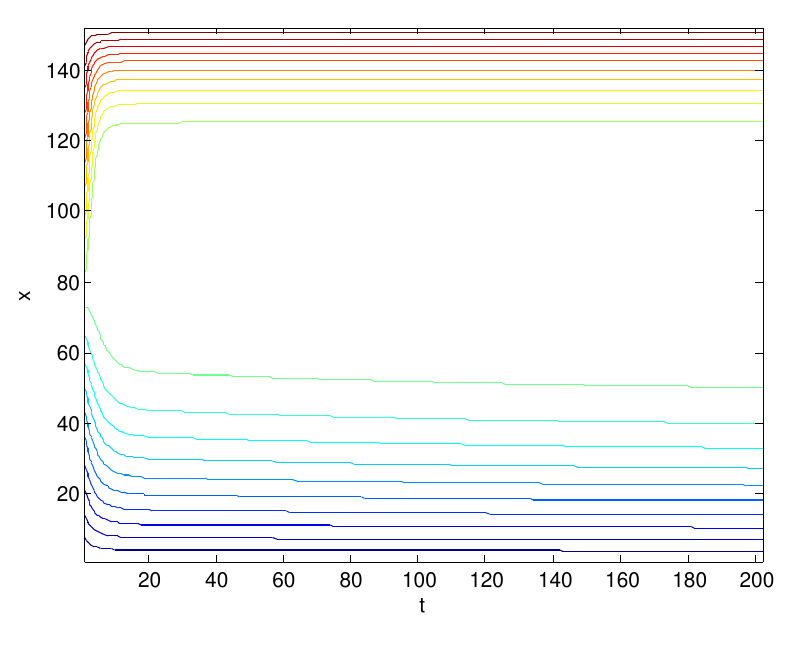}}
  \end{center}
\caption{\small{The left picture shows a numerical solution for $\lambda^+=3$,  $\lambda^-=1$ and $g(x)=16x-8.$ The right picture shows the contours of $u(t,x).$}}\label{1,3}
\end{figure}

  \begin{figure}[!htbp]
%{\includegraphics[width=7.0cm]{triangel}}
%\centering
 \begin{center}
  \subfloat[Numerical solution]{\includegraphics[width=.46\textwidth]{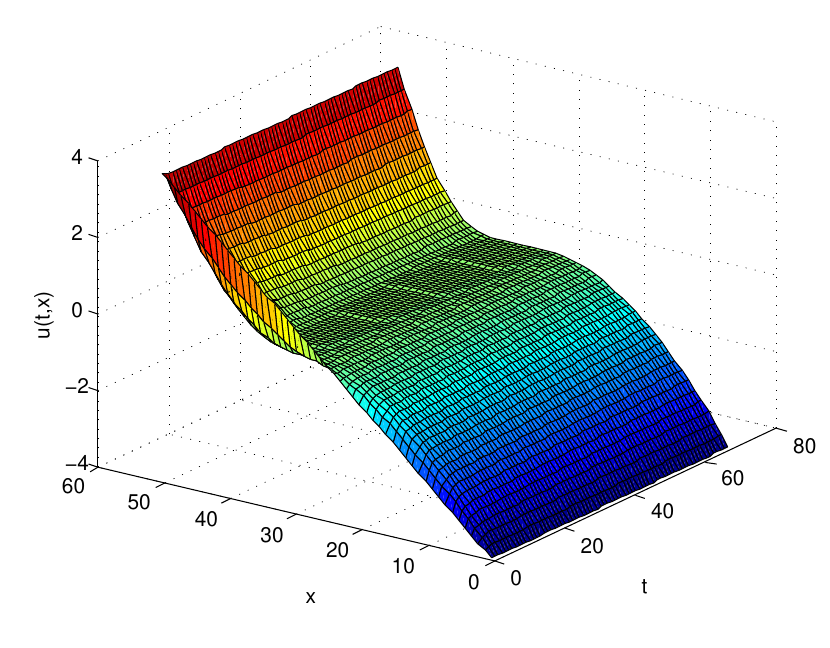}}
  \hspace{1cm}
  \subfloat[Contours]{\includegraphics[width=.4\textwidth]{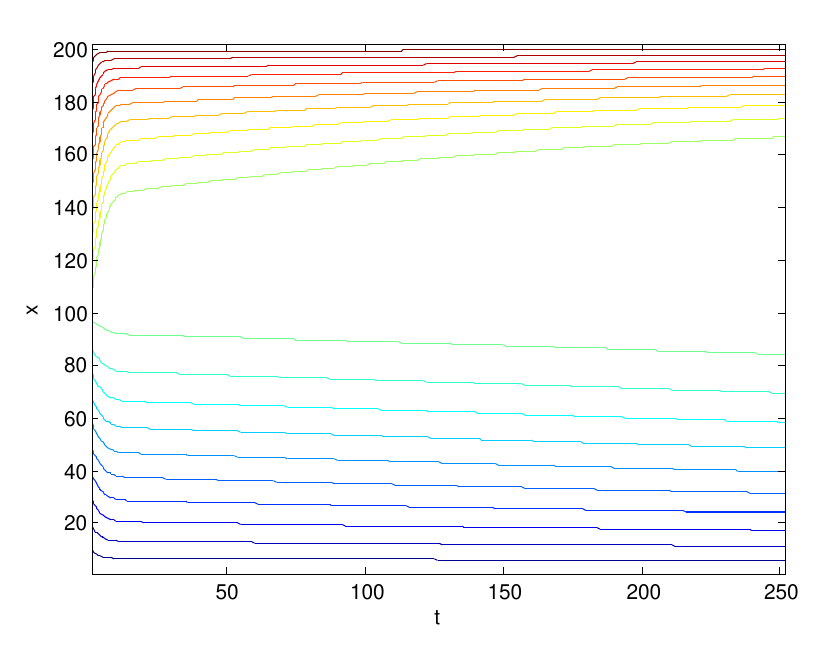}}
  \end{center}
\caption{\small{The left picture shows a numerical solution for $\lambda^+=0.7$,  $\lambda^-=0.2$ and $g(x)=8x-4.$ The right picture shows the contours of $u(t,x).$} }\label{7,2}
\end{figure}

  \begin{figure}[!htbp]
%{\includegraphics[width=7.0cm]{triangel}}
%\centering
 \begin{center}
  \subfloat[Numerical solution]{\includegraphics[width=.46\textwidth]{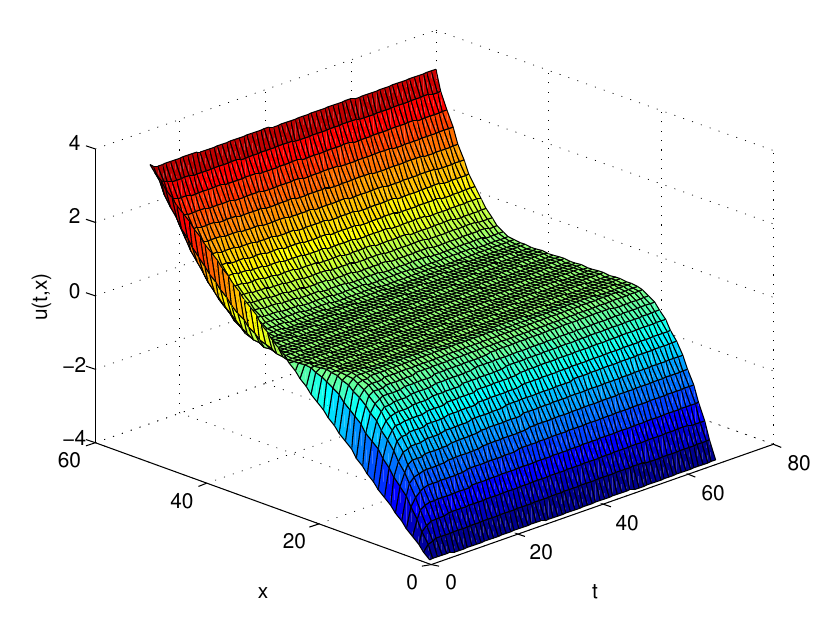}}
  \hspace{1cm}
  \subfloat[Contours]{\includegraphics[width=.4\textwidth]{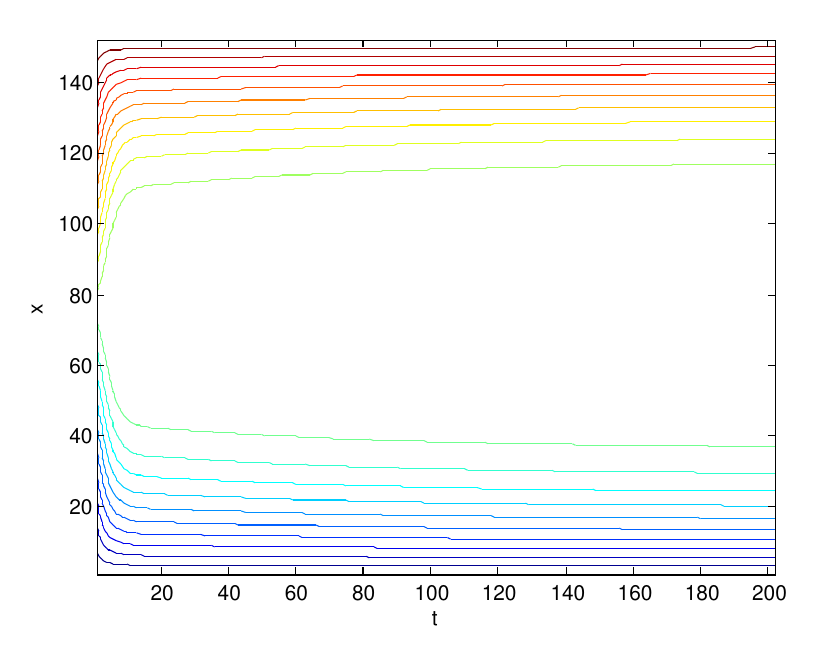}}
  \end{center}
\caption{\small{The left picture shows a numerical solution for $\lambda^+=\lambda^-=0.6$ and $g(x)=8x-4.$ The right picture shows the contours of $u(t,x).$} }\label{6,6}
\end{figure}

\section{Conclusion}
In this work we apply the theory of viscosity solutions to the parabolic two-phase obstacle-like problem in order to develop a convergent numerical scheme. Our developed schemes happened to be monotone, which allowed us to apply Barles-Souganidis theory to obtain their convergence. We observed that the method is applicable for the elliptic case as well. In section $4$  we applied the numerical method for different $\lambda^\pm,$ and concluded that the developed theory works well.

%\begin{figure}[!htbp]
%\includegraphics[width=.65\columnwidth]{lamda_2,7,16}
%\centering
%  \caption{Numerical solution of $u(t,x),$ for $g(x)=16(x-1/2)^3.$}\label{2,7,16}
%  \label{fig:E_2}
% \end{figure}

%\newpage
\bibliographystyle{acm}%

\bibliography{two_phase_parabolic}\label{sec-Ref}

\begin{thebibliography}{10}

\bibitem{MR1111111}
{\sc {Arakelyan}, A., {Barkhudaryan}, R., and {Poghosyan}, M.}
\newblock {Numerical Solution of the Two-Phase Obstacle Problem by Finite
  Difference Method}.
\newblock {\em ArXiv e-prints\/} (2014).

\bibitem{MR2961456}
{\sc Arakelyan, A.~G., Barkhudaryan, R.~H., and Poghosyan, M.~P.}
\newblock Finite difference scheme for two-phase obstacle problem.
\newblock {\em Dokl. Nats. Akad. Nauk Armen. 111}, 3 (2011), 224--231.

\bibitem{MR1315000}
{\sc Barles, G., Daher, C., and Romano, M.}
\newblock Convergence of numerical schemes for parabolic equations arising in
  finance theory.
\newblock {\em Math. Models Methods Appl. Sci. 5}, 1 (1995), 125--143.

\bibitem{MR1115933}
{\sc Barles, G., and Souganidis, P.~E.}
\newblock Convergence of approximation schemes for fully nonlinear second order
  equations.
\newblock {\em Asymptotic Anal. 4}, 3 (1991), 271--283.

\bibitem{MR1115934}
{\sc Bozorgnia, F.}
\newblock Numerical solutions of a two-phase membrane problem.
\newblock {\em Applied Numeric. Math. 61\/} (2011), 92--107.

\bibitem{MR1118699}
{\sc Crandall, M.~G., Ishii, H., and Lions, P.-L.}
\newblock User's guide to viscosity solutions of second order partial
  differential equations.
\newblock {\em Bull. Amer. Math. Soc. (N.S.) 27}, 1 (1992), 1--67.

\bibitem{MR0521262}
{\sc Duvaut, G., and Lions, J.-L.}
\newblock {\em Inequalities in mechanics and physics}.
\newblock Springer-Verlag, Berlin, 1976.
\newblock Translated from the French by C. W. John, Grundlehren der
  Mathematischen Wissenschaften, 219.

\bibitem{MR2532350}
{\sc Hu, B., Liang, J., and Jiang, L.}
\newblock Optimal convergence rate of the explicit finite difference scheme for
  {A}merican option valuation.
\newblock {\em J. Comput. Appl. Math. 230}, 2 (2009), 583--599.

\bibitem{MR1341739}
{\sc Ishii, H.}
\newblock On the equivalence of two notions of weak solutions, viscosity
  solutions and distribution solutions.
\newblock {\em Funkcial. Ekvac. 38}, 1 (1995), 101--120.

\bibitem{MR2218974}
{\sc Oberman, A.~M.}
\newblock Convergent difference schemes for degenerate elliptic and parabolic
  equations: {H}amilton-{J}acobi equations and free boundary problems.
\newblock {\em SIAM J. Numer. Anal. 44}, 2 (2006), 879--895 (electronic).

\bibitem{MR1934623}
{\sc Shahgholian, H.}
\newblock {$C\sp {1,1}$} regularity in semilinear elliptic problems.
\newblock {\em Comm. Pure Appl. Math. 56}, 2 (2003), 278--281.

\bibitem{MR2340105}
{\sc Shahgholian, H., Uraltseva, N., and Weiss, G.~S.}
\newblock The two-phase membrane problem---regularity of the free boundaries in
  higher dimensions.
\newblock {\em Int. Math. Res. Not. IMRN}, 8 (2007), Art. ID rnm026, 16.

\bibitem{MR2511041}
{\sc Shahgholian, H., Uraltseva, N., and Weiss, G.~S.}
\newblock A parabolic two-phase obstacle-like equation.
\newblock {\em Adv. Math. 221}, 3 (2009), 861--881.

\bibitem{MR2258264}
{\sc Shahgholian, H., and Weiss, G.~S.}
\newblock The two-phase membrane problem---an intersection-comparison approach
  to the regularity at branch points.
\newblock {\em Adv. Math. 205}, 2 (2006), 487--503.

\bibitem{MR1906034}
{\sc Uraltseva, N.~N.}
\newblock Two-phase obstacle problem.
\newblock {\em J. Math. Sci. (New York) 106}, 3 (2001), 3073--3077.
\newblock Function theory and phase transitions.

\bibitem{MR1825655}
{\sc Weiss, G.~S.}
\newblock An obstacle-problem-like equation with two phases: pointwise
  regularity of the solution and an estimate of the {H}ausdorff dimension of
  the free boundary.
\newblock {\em Interfaces Free Bound. 3}, 2 (2001), 121--128.

\end{thebibliography}

\end{document}